\documentclass[12pt,twoside,draft]{article}
\usepackage{amsmath,amssymb,mathrsfs,graphicx,bbm}

\usepackage{color}

\newtheorem{theorem}{Theorem}[section]
\newtheorem{lemma}[theorem]{Lemma}
\newtheorem{proposition}[theorem]{Proposition}
\newtheorem{corollary}[theorem]{Corollary}
\newtheorem{definition}[theorem]{Definition}

\newtheorem{rmrk}[theorem]{Remark}

\DeclareMathAlphabet{\mathbfit}{OML}{cmm}{b}{it}

\makeatletter
\@addtoreset{equation}{section}
\makeatother

\newenvironment{remark}
{\begin{rmrk} \em}
{\end{rmrk}}

\newcommand{\fn} {function}
\newcommand{\me} {measure}
\newcommand{\tr} {trajector}
\newcommand{\erg} {ergodic}
\newcommand{\sy} {system}

\newcommand{\R} {\mathbb{R}}

\newcommand{\Z} {\mathbb{Z}}
\newcommand{\N} {\mathbb{N}}
\newcommand{\qed} {\hfill {\small Q.E.D.} \par\medskip}
\newcommand{\skippar} {\par\medskip}
\newcommand{\ds} {\displaystyle}
\newcommand{\proof} {\noindent \textsc{Proof.} }
\newcommand{\proofof}[1] {\noindent \textsc{Proof of {#1}.} }

\newcommand{\eps} {\varepsilon}

\newcommand{\rset}[2] {\left\{ #1 \: \left| \: #2 \right. \! \right\} }
\newcommand{\lset}[2] {\left\{ \left. \! #1 \: \right| \: #2 \right\} }

\newcommand{\symmdiff} {\triangle}
\newcommand{\into} {\longrightarrow}

\renewcommand{\emptyset} {\varnothing}

\newcommand{\pr} {probability}
\newcommand{\ra} {random}
\newcommand{\rw} {random walk}
\newcommand{\en} {environment}
\newcommand{\sgn} {\mathrm{sgn}}
\newcommand{\dsy} {dynamical system}
\newcommand{\mj} {M}   
\newcommand{\Oen} {\Omega_\mathrm{en}}   
\newcommand{\nb} {\bar{n}}   

 \newcommand{\be}{\begin{equation}}
 \newcommand{\ee}{\end{equation}}

 \newcommand{\ba}{\begin{array}}
 \newcommand{\ea}{\end{array}}

 \newcommand{\bea}{\begin{eqnarray}}
 \newcommand{\eea}{\end{eqnarray}}

 \newcommand{\bl}{\begin{lemma}}
 \newcommand{\el}{\end{lemma}}

 \newcommand{\br}{\begin{remark}}
 \newcommand{\er}{\end{remark}}

 \newcommand{\bt}{\begin{theorem}}
 \newcommand{\et}{\end{theorem}}

 \newcommand{\bd}{\begin{definition}}
 \newcommand{\ed}{\end{definition}}

 \newcommand{\bcl}{\begin{claim}}
 \newcommand{\ecl}{\end{claim}}

 \newcommand{\bp}{\begin{proposition}}
 \newcommand{\ep}{\end{proposition}}

 \newcommand{\bc}{\begin{corollary}}
 \newcommand{\ec}{\end{corollary}}

 \newcommand{\bpr}{\begin{proof}}
 \newcommand{\epr}{\end{proof}}

 \newcommand{\bi}{\begin{itemize}}
 \newcommand{\ei}{\end{itemize}}

 \newcommand{\ben}{\begin{enumerate}}
 \newcommand{\een}{\end{enumerate}}

 \def \P {{\mathbb P}}
 \def \E {{\mathbb E}}

 \def \cA {\mathcal{A}}
 \def \cB {\mathcal{B}}
 \def \cC {\mathcal{C}}

 \def \cN {\mathcal{N}}

 \def \a {{\alpha}}

 \def \e {{\varepsilon}}

 \def \s {{\sigma}}
 
 \def \z {{\z}}
 \def \m {{\mu}}
 
 \def \t {{\tau}}

 \def \o {{\omega}}

 \def \x {{\xi}}
 \def \z {{\zeta}}

\def \à {{\`{a}}}
\def \ì{{\`{\i}}}
\def \ò{{\`{o}}}
\def \è{{\`{e}}}
\def \ù{{\`{u}}}

 \def \1{\mathbbm{1}} 

\begin{document}

\title{\textbf{Random walks in a one-dimensional \\
L\'evy random environment}}

\author{
\scshape
Alessandra Bianchi\,\thanks{
Dipartimento di Matematica Pura e Applicata,
Universit\`a di Padova, Via Trieste 63,
35121 Padova, Italy. E-mail:
\texttt{bianchi@math.unipd.it}.}
,
Giampaolo Cristadoro\,\thanks{
Dipartimento di Matematica, Universit\`a di Bologna,
Piazza di Porta San Donato 5, 40126 Bologna, Italy.
E-mails: \texttt{giampaolo.cristadoro@unibo.it,
marco.lenci@unibo.it}.}
, \\ \scshape
Marco Lenci\,$^\dagger$\thanks{
Istituto Nazionale di Fisica Nucleare,
Sezione di Bologna, Via Irnerio 46,
40126 Bologna, Italy.}
,
Marilena Ligab\`o\,\thanks{
Dipartimento di Meccanica, Matematica e Management,
Politecnico di Bari,
Via Orabona 4, 70125 Bari, Italy. E-mail:
\texttt{marilena.ligabo@poliba.it}.}
}

\date{Final version for the \\
  \emph{Journal of Statistical Physics} \\[4pt]
  January 2016}

\maketitle
\begin{abstract}
  We consider a generalization of a one-dimensional stochastic process
  known in the physical literature as L\'evy-Lorentz gas. The process
  describes the motion of a particle on the real line in the presence
  of a random array of marked points, whose nearest-neighbor distances
  are i.i.d.\ and long-tailed (with finite mean but possibly infinite variance).
  The motion is a continuous-time, constant-speed interpolation of a
  symmetric random
  walk on the marked points. We first study the quenched
  random walk on the point process, proving the CLT and the
  convergence of all the accordingly rescaled moments.
  Then we derive the quenched and annealed CLTs for the
  continuous-time process.

  \bigskip\noindent
  Mathematics Subject Classification (2010): 60G50, 60F05 (82C41, 60G55).
\end{abstract}

\section{Introduction}

For as long as they have existed, \rw s  have been used as models for
a wide range of transport processes in fields as
diverse as physics, chemistry and biology.

For a homogeneous \rw\  on a lattice,
under the hypothesis of finite variance of the distribution of jumps, classical
results include the central limit theorem (CLT), the functional CLT (a.k.a.\
invariance principle), and normal diffusion, defined as an asymptotically
linear time-dependence of the variance of the walker's position.

While the success of homogeneous random walks in capturing the main
features of transport in regular media is nowadays apparent, in many
interesting situations the walker moves in a complex and/or disordered
environment.
In such cases, correlations induced by spatial inhomogeneities can have a strong
impact on the transport properties, which cannot be simulated by a simple
homogenous model \cite{G, vB, vBS}.
This led (already 40 years ago) to the definition of a class of processes called
random walks in random environment (RWRE), where the transition probabilities
are themselves \ra\ \fn s of space (cf.\ \cite{Zei} for a review). This rich
class of walks is typically studied from two different viewpoints: that of the
\emph{quenched} processes, where one focuses on the dynamics for a typical fixed
\en, and that of the \emph{annealed} (or \emph{averaged}) processes,
where the interest is on the
effect of averaging over the environments.

On a related note, recent years have witnessed a growing interest around
anomalous diffusive processes, where the variance of a moving
particle has a super- or sub-linear growth in time. In the physical literature,
such anomalous behavior has been observed in many systems:
Lorentz gases with infinite horizon, rotating flows, intermittent dynamical
systems, etc.\ \cite{VHC}.

Several models
have been put forth
to describe such situations. Undoubtedly, the simplest among them are the
homogeneous \rw s whose transition probabilities have
an infinite second moment (and possibly an infinite first moment too)
\cite{GK}. Especially in the physical literature, they are sometimes dubbed
\emph{L\'evy flights}. Though
L\'evy flights
easily break normal diffusion, their defining feature is also their most serious
drawback, in that the variance of the walker's position is infinite at
\emph{all} times, failing to
reproduce the superlinear time-dependence that is typical of many systems
of interest, such as those mentioned earlier. More realistic models are
then considered, called \emph{L\'evy walks}: here the jumps are still picked from a
long-tailed distribution but the walker needs a certain time to complete a jump
(typically a time proportional to the length of the jump, implying
constant speed) \cite{CGLS, ZDK}.

Not much work has been done on systems that combine long-tailed
jumps and disordered media. To the authors' knowledge, the first
such examples are the L\'evy flights perturbed
by random drift fields introduced in \cite{F}. In this case
the cause of the anomalous diffusion is the
distribution of the jumps. Two more recent models are those of
\cite{BFK} and \cite{S}; though rather different from one another, both
systems are defined by a ``normal'' (meaning, simple, standard)
dynamics on an ``anomalous''
\en, which forces long jumps and is therefore responsible for the anomalous
behavior. In this sense, the models are representative of the many
physical situations (human mobility, epidemics, network routing, etc.)\
in which anomalous diffusion is caused by the complexity of
an underlying network (such as a small-world network).
The system presented in \cite{BFK}, called by the authors
\emph{L\'evy-Lorentz gas}, is the starting point of our investigation;
we will come back to it momentarily.
The only examples of long-tailed \rw s in \ra\ \en\
these authors have found in the rigorous mathematical
literature are the long-range walks on point processes studied
in \cite{BeR,CFG,Rou}.

A surge of interest in this topic has lately come from the physics of materials,
since a new glassy material has been devised, through which light
exhibits anomalous properties that can be experimentally controlled \cite{BB}.
The design of this so-called \emph{L\'evy glass} suggests an interpretation
of the motion of light in it by way of a L\'evy walk in a disordered environment,
as studied in \cite{ABG, ABG2, BCV, BBV, BDLV} (with varying
degrees of approximation). These papers focus
on the annealed versions of the models, and no rigorous proof is given.

Inspired by the above models, the system we study here is a generalization
of the L\'evy-Lorentz gas mentioned earlier. A \ra\ array of points,
called \emph{targets}, is given on the real line, such that
the distances between two neighboring targets are i.i.d.\ with finite
mean; they are, however, allowed to have infinite variance, which is the
interesting case here. A particle moves with unit speed between the
targets, driven by a \rw\ that is independent of all the rest.
More in detail, we assume that the origin is always a target
and that the particle starts from there. A \ra\ integer $\x$ is drawn
from a given distribution, upon which the particle starts to travel towards the
$\x^\mathrm{th}$ target. When the target has been reached the
procedure repeats from there.

This is therefore a continuous-time \rw, whose \tr ies have long inertial
segments due to a \ra\ \en, which is why we speak of
a \emph{\rw\ in a L\'evy \ra\ \en}. These are our results: for the (discrete-time)
\rw\ on the point process,
 we prove the quenched CLT and the convergence of all the normally
rescaled moments to those of a suitable Gaussian.
Then, by comparison, we derive  the quenched CLT for the
continuous-time process.
These results imply the annealed CLT for both the continuous- and
discrete-time walks.

The paper is organized as follows. In Section \ref{subs-model} we give
the precise definitions of all the processes associated with our \rw.
In Section \ref{subs-results}
we present our main results, whose proofs are found in Section
\ref{sec-proofs}, although some technical results are gathered in
Appendix \ref{app-tech}. Section \ref{subs-pvp} presents a construction
that is also of independent interest: a \dsy\ describing
the annealed process from the point of view of the particle.

\section{Model and main results}
\label{sec-mod-res}

\subsection{Definition of the model}
\label{subs-model}

We start by defining the following marked
point process on $\R$: let $\z := (\z_j, \, j\in\Z )$ be
a sequence of i.i.d.\ positive random variables with finite mean $\mu$,
and define the variables $\o_k$, $k\in\Z$, via
\be\label{PP}
  \o_0 := 0 \, , \qquad \o_k := \o_{k-1}+\z_{k} \,.
\ee
The process $\o := ( \o_k, \, k\in\Z )$ will be also referred to as
the \emph{\en}, and the single points $\o_k$ as the \emph{targets}.
We denote the set of all possible \en s by $\Oen$,
and the law just defined on it by $P$.

We are particularly interested in long-tailed $\z_j$, with infinite
variance, distributed for instance in the basin of attraction of
an $\alpha$-stable distribution, with $\alpha \in (1,2)$.
Environments of this type are usually called \emph{L\'evy \en s} in the
physical literature \cite{BFK, BCV, ABG}.

In order to define our continuous-time process, we need to
introduce two intermediate \rw s (RWs).
Let $ \Z^{\pm}$  be the positive/negative integers,
and $\N$ the non-negative integers.
Take $\x := (\x_i, \, i\in \Z^+)$ to be a sequence of i.i.d.\ $\Z$-valued
random variables, with density $p := (p_k, \, k \in \Z )$, where $p_k =
p_{-k}$ (symmetry condition), $p_{k+1} \le p_k$, $\forall k \ge 0$
(half-monotonicity), and such that
\be \label{vp}
  v_p := \sum_k k^2 p_k \in (0,\infty) \,.
\ee
Denote by $S := (S_n,\, n \in \N)$ the RW with increments provided by
the $\x_i$, that is
\be \label{underRW}
  S_0 := 0\,, \qquad  S_n := \sum_{i=1}^n \x_i \,, \quad \mbox{for } n
  \ge 1 .
\ee
This is called the \emph{underlying random walk}. It is
defined on the probability space $(\Z^\N, Q)$, endowed with the
$\s$-algebra generated by cylinder functions.

The second RW is defined, for each environment $\o \in \Oen$, as
\be\label{Y}
  Y_n \equiv Y_n^\o := \o_{S_n}\,, \quad \mbox{for } n \in \N \,.
\ee
In rough terms, $Y := (Y_n,\, n \in \N)$ performs the same jumps as $S$,
but on the points of $\o$. We call it the \emph{\rw\ on the point process}.
The associated \pr\ space is $(\o^\N, Q_\o)$, where $Q_\o$ is
the \pr\ induced on $\o^\N$ by $Q$ via (\ref{Y}) (more precisely,
$Q_\o$ is defined on the $\s$-algebra generated by cylinder \fn s).
In particular, for all $n\in\N$ and $k \in \Z$,
\be\label{measures}
  Q_\o(Y_n = \o_k) = Q(S_n =k) \,.
\ee

Once we fix the \en\ and the realization of the dynamics, that is,
for any given pair $(S, \o) \in (\Z^\N, \Oen)$, we can define the
sequence of \emph{collision times}
$\t(n) \equiv \t(n; S, \o)$ via
\be\label{collisiontime}
  \t(0) := 0 \,, \qquad \t(n):= \sum_{k=1}^n
  |\o_{S_k}-\o_{S_{k-1}}| \,, \quad \mbox{for } n\ge 1 \, .
\ee
Notice that, since the length of the $n^\mathrm{th}$ jump of the
walk is given by $|\o_{S_n}-\o_{S_{n-1}}|$,
$\t(n)$ represents the global length of the \tr y up to time $n$.

Finally, the process we are interested in is the continuous-time process
$X(t) \equiv X^\o(t)$ defined by
\be\label{process}
  X(t) := Y_n+ \sgn(\x_{n+1} )(t-\t(n)) \,, \quad \mbox{for } t \in
  [\t(n),\t(n+1)) \,,
\ee
where $\sgn$ is the sign \fn. In other words, $X:=( X(t) \,, t \in [0,\infty) )$
is the process whose trajectories interpolate those of the walk $Y$
and whose speed is 1 (save at collision times). In light of the discussion
made in the introduction, we describe the above as
a continuous-time \rw\ on a L\'evy \ra\ \en.

\begin{remark} \label{rk-lazy}
  Notice that $\x_{n+1} = 0 \Leftrightarrow S_{n+1} = S_{n} \Leftrightarrow
  \t(n+1) = \t(n)$. Therefore (\ref{process}) is never used in the case
  $\x_{n+1} = 0$, which makes the definition of $\sgn(0)$
  irrelevant there. More importantly, the self-jumps of the underlying RW
  (namely, $S_{n+1} = S_n$) are simply not seen by the process
  $X(t)$. This implies that we can remove any lazy component of
  $S = (S_n)$ by redefining
  \be\label{removelazy}
    p'_0 := 0 \,, \qquad p'_j := \frac{p_j} {\sum_{k \ne 0} p_k}, \
  \mbox{ for } j \ne 0.
  \ee
  (Notice that $\sum_{k \ne 0} p_k > 0$, because $v_p > 0$.)
  In particular, the case where $S$ is a simple symmetric RW,
  called L\'evy-Lorentz gas in \cite{BFK}, is included in our results.
\end{remark}

Indicate with $\cC := C([0,\infty); \R)$ the space of all continuous paths
from $[0,\infty)$ to $\R$, endowed with the Skorokhod topology.
We denote by $P_\o$ the \emph{quenched law} of $X$, which
is the probability induced by $Q$ on $\cC$ by the definitions
(\ref{Y})-(\ref{process}).

Finally, we use $\P$ for the \emph{annealed law} of the process,
defined on the space $\cC \times \Oen$ by
\be
\P(G \times F)=\int_F P_\o(G)P(d\o)\, .
\ee
This is the law that describes the entire randomness of the system.

\subsection{Main results}
\label{subs-results}

In order to state our main results we need to name a few parameters
pertaining to the underlying RW $S$. Let
\be\label{mj}
  \mj := \sum_{k \in \Z} |k| \, p_k = 2 \sum_{k = 1}^\infty k \, p_k
\ee
be its mean absolute jump, and denote
\be\label{qbar}
  \bar{q} := \sup \rset{ q \ge 0} {\sum_k |k|^q \, p_k < \infty} \,.
\ee
By our initial assumptions, $\bar{q} \ge 2$. The following
is standard:

\begin{proposition}\label{prop-standard}
  $S$ verifies the standard CLT, namely,
  \be\label{CLT-sn}
    \lim_{n\to\infty}\frac{S_n}{\sqrt n}\stackrel{d} = \cN(0, v_p).
  \ee
  Also, denoting by $E_Q$ the expectation w.r.t.\ $Q$ (the law of
  $S$),
  \be\label{rescmoment}
    m_q := \lim_{n \to \infty} \frac{E_Q ( |S_n|^q )}{n^{q/2}}
  \ee
  exists at least for all $q \ge 0$, $q \ne \bar{q}$. For $q \in [0, \bar{q})$,
  it is finite and equals the $q^{th}$ absolute moment of
  $\cN(0, v_p)$; for $q \in (\bar{q}, \infty)$, it is infinite.
\end{proposition}

For the sake of completeness,
we give the proof of Proposition \ref{prop-standard} at the beginning
of the next section. Observe that, if $\bar{q} > 2$, the  proposition says that
$(S_n / \sqrt{n})$ converges weakly to
a suitable Gaussian, together with all the moments of order $\le 2$.
If $\bar{q} = 2$, the proposition does not ensure that the second moment
converges, but its proof guarantees that $E_Q ( |S_n|^2 ) /n$ is
bounded above and below (this follows from (\ref{burkholder})-(\ref{st-10})
below and the fact that $E_Q ( |\x_1|^2 ) = v_p < \infty$ by hypothesis).
We describe this situation by saying that the underlying RW is
\emph{totally diffusive}.

\skippar

The purpose of this paper is to show that the \rw\ on the point process
$Y$ is also totally diffusive and its continuous-time interpolation $X$
verifies the CLT, both in the quenched sense, i.e., in a fixed
\en, for almost every \en. Recalling that $\mu$ denotes the mean of the
random variables $\zeta_i$,  these are our results:

\begin{theorem} \label{prop:qCLTY}
  For $P$-a.e.\ $\o\in\Oen$,
  \be\label{qCLTY}
    \lim_{n\to\infty} \frac{Y_n}{\sqrt{n}}\stackrel{d} = \cN(0, \mu^2\, v_p) \,.
  \ee
  The convergence is in distribution, relative to the law $P_\o$ on
  $\cC$.
\end{theorem}

\begin{theorem}\label{prop:qmomentY}
  Let $E_\o$ denote the expectation w.r.t.\ the \me\ $Q_\o$.
  If $q \in [0, \bar{q})$, then
  \be\label{pr:mom1}
    \lim_{n \to \infty} \frac{E_\o ( |Y_n|^q )} {n^{q/2}} = \mu^q m_q \,.
  \ee
  If also $q \in 2\N+1$, then
  \be\label{pr:mom2}
    \lim_{n \to \infty} \frac{E_\o ( Y_n^q )} {n^{q/2}} = 0 \,.
  \ee
  Both statements hold for $P$-a.e.\  $\o\in\Oen$.
\end{theorem}

\begin{theorem}\label{Th:qCLT}
  For $P$-a.e.\ $\o \in \Oen$,
  \be\label{qCLT}
    \lim_{t\to\infty}\frac{X(t)}{\sqrt{t}} \stackrel{d} = \cN \! \left
    (0, \frac{\mu}{\mj} v_p \right) \,.
  \ee
  The convergence is in distribution, relative to the law $P_\o$ on
  $\cC$.
\end{theorem}

\begin{remark} \label{rk-lazy2}
  In view of Remark \ref{rk-lazy}, let us observe that
  Theorem \ref{Th:qCLT} must not depend of the choice
  of $p_0$, the lazy component of $S$. However, as per
  definitions (\ref{vp}) and (\ref{mj}), $v_p$ and $\mj$
  do. On the other hand, if we define $v_p'$ and $\mj'$
  by using $(p_k')$ in lieu of $(p_k)$ in (\ref{vp}), (\ref{mj}),
  it is immediate to check that $v_p' / \mj' = v_p / \mj$.
\end{remark}

The quenched CLTs easily imply the annealed CLTs:

\begin{corollary} \label{co:ann-clt}
  The limits (\ref{qCLTY}) and (\ref{qCLT}) hold for the annealed processes
  as well, that is, w.r.t.\ the \me\ $\P$ on $\cC \times \Oen$.
\end{corollary}

Theorems \ref{prop:qCLTY} and \ref{prop:qmomentY}
provide a complete characterization of the quenched process $Y$.
As for the physically more relevant process $X$, it is an open question
whether a similar scaling for the quenched moments holds.
In the annealed case, heuristic arguments and numerical simulations
suggest that
the second moment does not always grow linearly in time. In particular
\cite{BCV},
if the distribution of the distance between targets behaves
like $dP (\z_0 \le z) / dz \sim z^{-1-\a}$, for $z \to \infty$, the second
moment is expect to scale like
\be
  \E\! \left( X(t)^2 \right) \sim \left\{
  \begin{array}{ll}
    t^{5/2 - \a} \,, & 1 \le \a \le 3/2 \,; \\
    t \,, & \a > 3/2 \,.
  \end{array} \right.
\ee

\section{Proofs}
\label{sec-proofs}

In this section we prove our results. The most elaborate proof,
that of Theorem \ref{Th:qCLT}, requires a representation of the annealed
process
as a \dsy\ `from the point of view of the particle'. We present
this \sy\ in Section \ref{subs-pvp}. The technical lemmas which
offer little insight on the flow of the proofs are given in the
Appendix.

We start with the standard results about the underlying RW.
\bigskip

\proofof{Proposition \ref{prop-standard}} The CLT is a well-known
result for a finite-variance RW. The convergence of the
rescaled moments is in general not so well-known. For this reason
we provide a short proof, though other references
may be found in the literature (e.g., \cite{Y}).

Denoting $\chi_n := \left( \sum_{i=1}^n \x_i^2 \right)^{1/2}$,
Burkholder's inequality \cite[Thm.~3.2]{B} states that, for all
$q>1$, there exist constants $C_q > c_q > 0$ such that
\be\label{burkholder}
  c_q \| \chi_n \|_q \le \| S_n \|_q \le C_q \| \chi_n \|_q \,,
\ee
where $\| \cdot \|_q$ denotes the $L^q$-norm w.r.t.\ $Q$.
For $1 < q < \bar{q}$,
$E_Q (  |\x_i|^q ) := \sum_k |k|^q \, p_k < \infty$. By (\ref{burkholder}),
$E_Q ( |S_n|^q ) = \| S_n \|_q^q$ is
asymptotic to
\be\label{st-10}
  \left\| \chi_n \right\|_q^q = \left\| \sum_{i=1}^n \x_i^2
  \right\|_{q/2}^{q/2} \le \left( \sum_{i=1}^n \left\| \x_i^2
  \right\|_{q/2} \right)^{q/2} \!\! = \, n^{q/2} E_Q ( |\x_1|^q ) \,.
\ee
This implies that the $q^\mathrm{th}$
absolute moment of $(S_n / \sqrt{n})$ is bounded above in $n$.
Since the process converges weakly to $\cN(0, v_p)$, a
simple argument \cite[Ex.~2.5]{D}
shows that, $\forall q' \in [0,q)$, the $(q')^\mathrm{th}$ absolute moment
of $(S_n / \sqrt{n})$
converges to the $(q')^\mathrm{th}$ absolute moment of $\cN(0, v_p)$.
Since $q$ was arbitrary, the conclusion holds for all $q \in [0, \bar{q})$.

On the other hand, when $q \in (\bar{q}, \infty)$, $E_Q ( |S_n|^q )$
behaves like
\be\label{st-50}
   \left\| \sum_{i=1}^n \x_i^2 \right\|_{q/2}^{q/2} \ge
   E_Q \! \left( \sum_{i=1}^n |\x_i|^q \right) = \infty \,,
\ee
having used the inequality $(a+b)^{q/2} \ge a^{q/2} + b^{q/2}$ (holding for
$a,b \ge 0$ and $q/2 \ge 1$) and the fact that $E_Q ( |\x_i|^q ) = \infty$.
\qed

\subsection{The \rw\ on the point process}

\proofof{Theorem \ref{prop:qCLTY}} This proof follows that
of Thm.~1.13 of \cite{BeR}. By the definition (\ref{PP}) of $\o$, we have
\be\label{omega-n}
  \o_n = \left\{
  \begin{array}{ll}
    \sum_{k=1}^{n} \z_k \,, & n>0 \,; \\
    0 \,,  & n=0 \,; \\
    -\sum_{k=n+1}^{0} \z_k \,, & n<0 \,.
  \end{array}
  \right.
\ee
The strong law of large numbers on $(\z_k)$ can be expressed as follows:
fixed $b \in \R$, for $P$-a.e.\ $\o \in \Oen$,
\be\label{llnPP}
  \lim_{j \to \infty} \frac{\o_{[bj]}} {j} = b\mu \,,
\ee
where $[r]$ denotes the integer part of $r \in \R$.

Since $Y_n= \o_{S_n}$, we get that, for $a \in \R$, $\e > 0$ and
$P$-a.e.\ $\o \in \Oen$,
\begin{align}
  \lim_{n\to\infty} Q_\o \! \left( \frac{Y_n}{\sqrt n} \le a \right)
  &\le \lim_{n\to\infty} \left[ Q\! \left( \frac{\o_{S_n}}{\sqrt n} \le a \,,
  \frac{S_n}{\sqrt n} > \frac{a}{\m} + \e \right) + Q\! \left( \frac{S_n}
  {\sqrt n} \le \frac{a}{\m} + \e \right) \right] \nonumber \\
  &\le \lim_{n\to\infty} \left[ Q\! \left( \frac{\o_{[(\frac{a}{\m} + \e) \sqrt n]}}
  {\sqrt n} \le a \right) + Q\! \left( \frac{S_n}{\sqrt n} \le \frac{a}{\m} + \e
  \right) \right] \nonumber \\
  &= \Phi\! \left( \frac{a}{\mu \sqrt{v_p}} + \e' \right) \,,
\end{align}
where $\e' := \eps/\sqrt{v_p}$, $\Phi$ is the distribution function of the
standard normal, and we have used (\ref{llnPP}) and (\ref{CLT-sn}).
Analogously,
\be
  \lim_{n\to\infty} Q_\o\! \left( \frac{Y_n}{\sqrt n} \le a \right) \ge
  \Phi\! \left( \frac{a}{\mu\sqrt{v_p}} - \e' \right) \,,
\ee
and altogether one gets the desired convergence.
\qed

\skippar

\proofof{Theorem \ref{prop:qmomentY}} The basic ingredients of the
proof are the convergence of the moments of $(S_n / \sqrt{n})$ and the
law of large numbers for $\o_n$, cf.\ (\ref{llnPP}).

By the latter, for every $\e > 0$ and $P$-a.e.\ $\o$, there exists $k_0
\equiv k_0(\e,\o)$ such that, for all $|k| \ge k_0$,
\be\label{errorcontrol}
  \left|\frac{\o_k}{k}-\mu \right| < \e \,.
\ee
In particular there exists $c \equiv c(\o) > 0$ such that $|\o_k| \le c |k|$,
for all $k \in \Z$. Let us fix a $P$-typical $\o$. Recalling that $Y_n= \o_{S_n}$,
we have:
\be\label{uno}
\begin{split}
  \frac{ E_\o (|Y_n|^q) } {n^{q/2}}
  &= E_Q\! \left( \frac{|\o_{S_n}|^q}{|S_n|^q} \, \frac{|S_n|^q}{n^{q/2}}\right) \\
  &= E_Q\! \left( 1_{ \{ |S_n| \ge k_0 \} } \frac{|\o_{S_n}|^q}{|S_n|^q} \,
  \frac{|S_n|^q}{n^{q/2}} \right) +
  E_Q\! \left( 1_{ \{ |S_n| < k_0 \} } \frac{|\o_{S_n}|^q}{|S_n|^q} \,
  \frac{|S_n|^q}{n^{q/2}} \right) \\
  &\le (\mu+\e)^q \, E_Q\! \left( \frac{|S_n|^q}  {n^{q/2}} \right) + c^q \,
  Q ( |S_n|<k_0 )^{1/r'} E_Q\! \left( \frac{|S_n|^{rq}} {n^{rq/2}}
  \right)^{1/r} ,
\end{split}
\ee
where in the last step we have used H\"older's inequality
with $r>1$ and $1/r' + 1/r = 1$. Now choose $r$ so that
$rq<\bar q$. By Proposition \ref{prop-standard},
\begin{align}
  & \lim_{n\to\infty} Q\left( |S_n| < k_0\right) = 0 \,; \\
  & \lim_{n\to\infty} E_Q\! \left( \frac{|S_n|^q} {n^{q/2}} \right) = m_q \,; \\
  & \lim_{n\to\infty} E_Q\! \left( \frac{|S_n|^{rq}} {n^{rq/2}} \right) =
  m_{rq} < \infty \,.
\end{align}
In conclusion,
\be\label{fbound1}
  \limsup_{n\to\infty} \, \frac{E_\o ( |Y_n|^q )} {n^{q/2}}
  \le (\mu+\e)^q m_q \,.
\ee

Similarly, for the lower bound,
\be\label{momentstima2}
\begin{split}
  E_Q\! &\left( \frac{|\o_{S_n}|^q}{|S_n|^q} \, \frac{|S_n|^q}{n^{q/2}} 
  \right) \ge (\mu-\e)^q \, E_Q\! \left( 1_{ \{ |S_n| \ge k_0 \} } 
  \frac{|S_n|^q} {n^{q/2}} \right) \\
  &\quad = (\mu-\e)^q \, E_Q\! \left( \frac{|S_n|^q}{ n^{q/2}} \right) - 
  (\mu-\e)^q \,E_Q\! \left( 1_{ \{ |S_n| < k_0 \} } \frac{|S_n|^q} {n^{q/2}} 
  \right) \\
  &\quad \ge  (\mu-\e)^q \, E_Q\! \left(\frac{|S_n|^q}{ n^{q/2}}\right) - 
  (\mu-\e)^q \, Q ( |S_n|<k_0 )^{1/r'} E_Q\! 
  \left( \frac{|S_n|^{rq}} {n^{rq/2}} \right)^{1/r} ,
  \end{split}
\ee
which, by the same arguments as above, gives
\be\label{fbound2}
  \liminf_{n\to\infty} \, \frac{E_\o ( |Y_n|^q )} {n^{q/2}}
  \ge (\mu-\e)^q m_q \,.
\ee
Assertion (\ref{pr:mom1}) follows from (\ref{fbound1}), (\ref{fbound2})
and the arbitrariness of $\e$. The limit (\ref{pr:mom2}) is proved in a
similar fashion upon rewriting
\be
  \frac{ E_\o ( Y_n^q ) } {n^{q/2}} =
  \frac{ E_\o\! \left( 1_{ \{ S_n \ge 0 \} } \, |Y_n|^q \right) } {n^{q/2}} -
  \frac{ E_\o\! \left( 1_{ \{ S_n < 0 \} } |\, Y_n|^q \right) } {n^{q/2}} \,.
\ee
\qed

\subsection{The point-of-view-of-the-particle \dsy}
\label{subs-pvp}

In this section we introduce a process, or rather a \dsy, that
describes the point of view of the particle for the RW on the point
process.

Keeping in mind the definitions and notation of Section \ref{subs-model},
let $(\Z^{\Z^+}, Q_o)$ denote the \pr\ space of the sequences
$\x = ( \x_i, \, i \in \Z^+)$, namely, $Q_o$ is the \pr\ for
$\Z^+$ copies of i.i.d.\ $\Z$-valued variables with density $p =
(p_k, \, k \in \Z)$.
Indicate with $\sigma_\x$ the left shift on this space. Evidently,
$\sigma_\x$ preserves $Q_o$ and is \erg.
This process is isomorphic to the RW $(S_n, \, n \in \N)$, by
construction of the latter. When conjugated with the natural
isomorphism, $\sigma_\x$ acts on $(\Z^\N, Q)$ as:
$(S_n, \, n \in \N) \mapsto ( S_{n+1} - S_1, \, n \in \N)$.

Further denote by $((\R^+)^\Z, P_o)$ the \pr\ space of the sequences
$\z := (\z_j ,\, j \in \Z)$, where $P_o$ is the
Bernoulli \me\ based on the variables $\z_j$ defined
in Section \ref{subs-model}. Indicate with $\sigma_\z$ the left shift on this space:
$\sigma_\z$ is an \erg\ automorphism of $((\R^+)^\Z, P_o)$.
Again, there is a natural isomorphism between $((\R^+)^\Z, P_o)$
and $(\Oen, P)$. Upon conjugation with it,
$\sigma_\z$ acts on $(\Oen, P)$ as:
$(\o_k ,\, k \in \Z) \mapsto ( \o_{k+1} - \o_1 ,\, k \in \Z)$. Also,
$\sigma_\z^{-1}$ acts as: $(\o_k ,\, k \in \Z) \mapsto ( \o_{k-1} - \o_{-1}
,\, k \in \Z)$.

Set $\Sigma := \Z^{\Z^+} \times (\R^+)^\Z$ and $\nu := Q_o \otimes P_o$,
and define $T: \Sigma \into \Sigma$ via
\be \label{T}
  T(\x, \z) := (\sigma_\x(\x),  \sigma_\z^{\x_1} (\z)) \,.
\ee
We think of $(\Sigma, \nu, T)$ as a \dsy. Let us call $\x$ the
\emph{dynamical variable} and $\z$
the \emph{environmental variable}, or simply the \emph{environment}.
Fix an initial condition $(\x, \z)$. The first component of the dynamical
variable, $\x_1$, determines the jump that the underlying RW is about to
make, namely $Y_1 = \o_{\x_1}$.
Applying $T$ translates the environment by the
quantity $-Y_1$ (corresponding to $|\x_1|$ discrete shifts in the
direction opposite to the jump), and shifts the dynamical variable,
so that the \sy\ is ready for the next jump (determined by $\x_2$)
under the pretense that $Y_1$ is the origin.

In other words, this \dsy\ describes the
annealed process from the point of view of the particle (PVP). This
is why we call it the \emph{PVP \dsy}.

\begin{theorem}
  \label{thm-pvp}
  $(\Sigma, \nu, T)$ is \me-preserving and \erg.
\end{theorem}

The proof of this theorem is found in Appendix \ref{app-pvp}.
The isomorphisms $\x \leftrightarrow S$
and $\z \leftrightarrow \o$, mentioned earlier,
entail that Theorem \ref{thm-pvp} is equivalent to the following:

\begin{corollary}
  \label{cor-pvp}
  The mapping
  \be\label{alt-pvp-map}
    \left( (S_n ,\, \o_k), \, n \in \N, \, k \in \Z \right) \mapsto
    \left( (S_{n+1} - S_1 ,\, \o_{k+S_1} - \o_{S_1}), \, n \in \N,
    \, k \in \Z \right)
  \ee
  on $(\Z^\N \times \Oen, Q \otimes P)$ is \me-preserving and \erg.
\end{corollary}

The following technical lemma, needed in the proof of
the main result, will also be proved in Appendix \ref{app-pvp}.

\begin{lemma}
  \label{lem-pvp-me}
  For $\z \in (\R^+)^\Z$, set $F_\z := \Z^{\Z^+} \times \{ \z \}$
  (in the remainder, any such set will be referred to as an
  horizontal fiber of $\Sigma$). Then,
  \be\label{pvp-me-1}
    T (F_\z) = \bigcup_{k \in \Z} F_{\sigma_\z^k (\z)} \,.
  \ee
  Furthermore (with a minor abuse of notation) indicate with
  $Q_o ( \, \cdot \,|\, F_\z)$ the \me\ on $F_\z$ induced by $Q_o$
  via the identification $F_\z \cong \Z^{\Z^+}$. $T$ pushes this
  \me\ to
  \be\label{pvp-me-2}
    T_* Q_o ( \, \cdot \,|\, F_\z) =
    \sum_{k \in \Z} p_k \, Q_o ( \,\cdot\, |  F_{\sigma_\z^k (\z)} ) \,.
  \ee
\end{lemma}

\subsection{CLT of the L\'evy walk}

We will prove Theorem \ref{Th:qCLT} by controlling the
continuous-time walk $(X(t))$ through the discrete-time
walk $(Y_n)$. To this goal, it is convenient to
introduce a quantity which counts
the number of collisions of the process $X(t)$ up to time $t$.
Formally, for every $t \in \R^+$, set
\be\label{collisionnumber}
  n(t) \equiv n(t; S, \o) := \max \rset{ m\in\N} { t\ge \t(m) } \,.
\ee
This is a sort of inverse \fn\ of the collision time $\t(n)$,
defined in (\ref{collisiontime}). In point of fact, when $\t(n)$
is strictly monotonic (which occurs when $S$ has no lazy component,
cf.\ Remark \ref{rk-lazy}), $n(t)$ is a suitable piecewise extension
of the inverse of $\t(n)$.

\begin{lemma} \label{lemma:time}
In view of the definitions (\ref{collisiontime}) and (\ref{collisionnumber}),
which depend on $(S, \o) \in \Z^\N \times \Oen$, we have that,
$(Q \otimes P)$-almost surely, equivalently, $\P$-almost surely,
\begin{align}
  \label{time1}
  \lim_{n \to \infty} & \frac{\t(n)}{n} = \mj \mu \,; \\
  \label{time2}
  \lim_{t \to \infty} \: & \frac{t}{n(t)} = \mj \mu \,.
\end{align}
\end{lemma}

\proof By (\ref{collisiontime}) we see that
$\t(n)$ is the Birkhoff sum of the \fn\
\be\label{time-10}
g(S,\o) := |\o_{S_1} - \o_{S_0}| = |\o_{S_1}| = |\o_{\x_1}| \,,
\ee
on $\Z^\N \times \Oen$, relative to the dynamics
(\ref{alt-pvp-map}). So (\ref{time1}) follows by Corollary \ref{cor-pvp}
and the Birkhoff theorem: for $(Q \otimes P)$-a.e.\ choice of
$(S, \o) \in \Z^\N \times \Oen$,
\be\label{time-20}
\begin{split}
  \lim_{n \to \infty} \frac{\t(n)}{n} &=
  \int_{\Z^\N \times \Oen} \!\! g \: d(Q \otimes P) \\
  &= \int_{\Z^{\Z^+} \times \Oen}
  \!\!\ |\o_{\x_1}| \: d(Q_o \otimes P) \\
  &= \sum_{k \in \Z} p_k \int_{\Oen} \!\!\ |\o_k| \: d P \\
  &= \sum_{k \in \Z} p_k |k| \mu = \mj \mu \,,
\end{split}
\ee
having used some of the notation and arguments given
in Section \ref{subs-pvp}.

Moreover, since by definition $n(t) \to \infty$, almost surely, as $t \to \infty$,
(\ref{time1}) implies that
\be
\lim_{t \to \infty} \frac{\t(n(t)+1)} {n(t)} = \lim_{t \to \infty}
\frac{\t(n(t))}{n(t)} = \mj \mu \,.
\ee
But
\be
  \lim_{t \to \infty} \left| \frac{t}{n(t)} - \frac{\t(n(t))}{n(t)} \right| \le
  \lim_{t \to \infty} \left( \frac{\t(n(t)+1)}{n(t)} - \frac{\t(n(t))}{n(t)} \right) = 0 \,,
\ee
giving (\ref{time2}).
\qed

\begin{lemma}\label{lem-mom}
  For $P$-a.e.\ $\o \in \Oen$ (equivalently, $P_o$-a.e.\ $\z \in (\R^+)^\Z$)
  and every $\ell \in \Z$,
  \be\label{mom-stat}
    \lim_{n \to \infty} E_\o \! \left( \o_{S_n + \ell} - \o_{S_n + \ell - 1}
    \right) = \lim_{n \to \infty} E_\o \! \left( \z_{S_n + \ell} \right) =
    \mu \,.
  \ee
  Also, given any positive sequence $\psi_n = o(\sqrt{n})$,
  $n \to \infty$, the above limits are uniform for $|\ell| \le \psi_n$
  (for a fixed $\o$, or $\z$).
\end{lemma}

\proof To start with, the first equality of (\ref{mom-stat}) follows
trivially by the definitions of $\o$ and $\z$, so we prove the second
one.

Again, we use the machinery of Section \ref{subs-pvp}. Set
$h (\x,\z) := \z_\ell$. By (\ref{T}) and (\ref{underRW}) we write
\be\label{mom-30}
  h \circ T^n (\x, \z) = h (\s_\x^n(\x), \s_\z^{S_n}(\z)) =
  \z_{S_n + \ell} \,.
\ee
Now, recall the definitions given in the statement of Lemma
\ref{lem-pvp-me}. Observe that taking the
expectation $E_\o$ is tantamount to integrating over the fiber
$F_\z$ w.r.t.\ the \me\ $Q_o ( \, \cdot \,|\, F_\z)$, where $\z$
corresponds to $\o$ via (\ref{PP}). Hence,
with the notation
 \be\label{p-n-j}
  p^{(n)}_j \!\!\!\ := \sum_{k_1 + \cdots + k_n = j} p_{k_1} \cdots  p_{k_n}
  = Q (S_n = j) \,,
\ee
we get
\be\label{mom-40}
\begin{split}
  E_\o \! \left( \z_{S_n + \ell} \right)
  &= \int_{F_\z}  (h \circ T^n) \, d Q_o ( \, \cdot \,|\, F_\z) \\
  &= \sum_{k_1, \ldots, k_n} p_{k_1} \cdots  p_{k_n} \int h \,
  d Q_o ( \, \cdot \,|\, F_{\sigma_\z^{k_1 + \cdots + k_n} (\z)} ) \\
  &= \sum_{j \in \Z} p^{(n)}_j \! \int h \, d Q_o ( \, \cdot \,|\,
  F_{\sigma_\z^j (\z)} ) \\
  &= \sum_{j \in \Z} p^{(n)}_j \, \z_{j + \ell} \,.
\end{split}
\ee
In the second equality above, we have applied Lemma
\ref{lem-pvp-me}
recursively $n$ times: the summation is over $\Z^n$ and each
integral is taken over the horizontal fiber specified by the
integration \me. In the fourth equality we have used that
$h$ is constant along horizontal fibers.

At this point we want to apply Lemma \ref{lem-tech1} of the Appendix
with $a_j := \z_{j + \ell}$
and $p^{(n)}$ as above. We need to check the hypotheses of the lemma.
First off, $p^{(n)}$ verifies condition \emph{(i)} because it is
symmetric and half-monotonic (this is, e.g., a consequence of
Lemma \ref{lem-tech2}, as $p$ is symmetric and half-monotonic
by assumption). It also verifies condition \emph{(ii)}, because
the underlying RW satisfies the CLT.

\begin{remark}
  This is the only point in the paper where the half-monotonicity of 
  $p$ is used.
\end{remark}

As for the hypothesis on $a$, we use the \erg ity of the process
$\z$. Thus, for $P_o$-a.e.\ $\z \in (\R^+)^\Z$,
\be\label{mom-45}
  \lim_{k \to \infty} \, \frac1 k  \sum_{j=0}^{k-1} a_j =
  \lim_{k \to \infty} \, \frac1 k \sum_{j=-1}^{-k} a_j =
  \E( \z_1 ) = \mu \,.
\ee
So Lemma \ref{lem-tech1}
can be applied almost surely in $\z$, equivalently in $\o$.
Using the notation of that lemma, (\ref{mom-40}) becomes
\be\label{mom-50}
  E_\o \! \left( \z_{S_n + \ell} \right) =
  \sum_{j \in \Z} p^{(n)}_j a_j = \mathcal{E}_n (a) \,,
\ee
showing that, for $P$-a.e.\ $\o \in \Oen$, (\ref{mom-50})
converges to  (\ref{mom-45}), as $n \to \infty$, which proves the
limit (\ref{mom-stat}).

Finally, observe that the underlying \rw\ is strongly aperiodic by 
hypothesis: this implies (rather easily) that, as $n \to \infty$, 
$p_{j-\ell}^{(n)} - p_j^{(n)} = o( p_j^{(n)} )$, uniformly
in $j \in \Z$ and $|\ell| \le \psi_n$. Since (\ref{mom-50}) can be rewritten
as $E_\o ( \z_{S_n + \ell} ) = \sum_j p^{(n)}_{j-\ell} \, \z_j$, its limit
is the same as for $\ell = 0$, uniformly for $|\ell| \le \psi_n$.
\qed

We are now ready to prove our main theorem (we will not
prove the obvious Corollary \ref{co:ann-clt}).
\bigskip

\proofof{Theorem \ref{Th:qCLT}}
Let us define $\nb(t) := [ t / \mj \mu ]$. Compare $\nb(t)$ to $n(t)$:
for fixed $t$, the former is a constant while the latter is a \ra\ variable
on $(\o^\N, Q_\o)$. For $P$-a.e.\ $\o \in \Oen$, we have that,
$Q_\o$-almost surely,
\begin{equation}\label{new-03}
  \lim_{t \to \infty} \frac{n(t) - \nb(t)} t = 0
\end{equation}
(this follows form Lemma \ref{lemma:time} and Fubini's Theorem).
Moreover, by Theorem \ref{prop:qCLTY} and the definition of
$\nb(t)$,
\be\label{new-05}
\lim_{t\to\infty} \frac{Y_{\nb(t)}}{\sqrt{t}} \stackrel{d} {=} \cN \! \left(0,
\frac{\mu}{\mj} v_p \right) \,,
\ee
for $P$-a.e.\ $\o$. Since $X(t)$ always lies between
$Y_{n(t)}$ and $Y_{n(t)+1}$, it is easy
to see that
\be\label{new-08}
\begin{split}
  \left| \frac{X(t)}{\sqrt{t}} - \frac{Y_{\nb(t)}}{\sqrt{t}} \right| &\le
  \max \left\{ \frac{| Y_{n(t)} - Y_{\nb(t)} |}{\sqrt{t}} \,, \frac{| Y_{n(t)+1} -
  Y_{\nb(t)} |}{\sqrt{t}} \right\} \\
  &\le \frac{| Y_{n(t)} - Y_{\nb(t)} |}{\sqrt{t}} + \frac{| Y_{n(t)+1} -
  Y_{\nb(t)} |}{\sqrt{t}} \,.
\end{split}
\ee
In light of (\ref{new-05}), and using Slutzky's Theorem
\cite[Thm.~13.18]{Kle},
Theorem \ref{Th:qCLT} will be proved once we prove that,
$P$-almost surely,
the two terms in the second line of (\ref{new-08}) converge
to 0 in distribution, and thus in \pr, w.r.t.\ $Q_\o$. We will only
show the convergence of the first term, the second one being
completely analogous.

Applying the Portemanteau Theorem \cite[Thm.~13.16]{Kle}, it will
suffice to prove that, given an $\o$ for which (\ref{new-03}) holds,
and a bounded Lipschitz \fn\ $f: \R \to \R$,
\be\label{new-20}
  \lim_{t\to \infty} E_\o \! \left( f \! \left(  \frac{|Y_{n(t)} - Y_{\nb(t)} |}
  {\sqrt{t}} \right) \right) = f(0) \,.
\ee
So, fix $\e > 0$. By (\ref{new-03}), one can find a `bad' set
$B_1 \subset \o^\N$, with $Q_\o(B_1) \le \e/6 \|f\|_\infty$, and
a \fn\ $\phi: \R^+ \into \R^+$, with
$\lim_{t \to \infty} \phi(t)/t = 0$, such that
\be
  |n(t) - \nb(t)| \le \phi(t)
\ee
for all realizations of the dynamics in $\o^\N \setminus
B_1$. Moreover, by (\ref{CLT-sn}), there exist another bad set
$B_2 \subset \o^\N$, again with $Q_\o(B_2) \le \e/6 \|f\|_\infty$,
and a constant $C > 0$ such that, for all
realizations in $\o^\N \setminus B_2$,
\be
  \left| S_{n(t)} - S_{\nb(t)} \right| \le C \sqrt{ | n(t) - \nb(t) |}.
\ee
Altogether, for all realizations in $\o^\N \setminus (B_1 \cup B_2)$,
\be\label{new-40}
  \left| S_{n(t)} - S_{\nb(t)} \right| \le C \sqrt{ \phi(t) }.
\ee

We split the average in the l.h.s.\ of (\ref{new-20}) in two parts, restricting
it, respectively, to $B_1 \cup B_2$ and its complement $G := \o^\N
\setminus (B_1 \cup B_2)$. For the first part, we estimate
\be\label{new-50}
  E_\o \! \left( 1_{B_1\cup B_2} \left| f \! \left(  \frac{|Y_{n(t)} - Y_{\nb(t)} |}
  {\sqrt{t}} \right) -f(0) \right| \right) \le 2 \|f\|_\infty \, Q_\o(B_1\cup B_2)
  \le \frac23 \e \,,
\ee
where $1_A$ denotes the indicator \fn\ of $A \subset \o^\N$.
For the second part, if $c$ is the Lipschitz constant of $f$, we
write
\be\label{new-55}
  E_\o \! \left( 1_G \left| f \! \left(  \frac{|Y_{n(t)} - Y_{\nb(t)} |}
  {\sqrt{t}} \right) -f(0) \right| \right) \le \frac{c}{\sqrt{t}} \, E_\o
  \! \left( 1_G \left| Y_{n(t)} - Y_{\nb(t)} \right| \right) \,.
\ee

By definition of the processes $Y$, $\o$, and $\z$ (cf.\
Section \ref{subs-model}),
\be\label{new-57}
  Y_{n(t)} - Y_{\nb(t)} = \o_{S_{n(t)}} - \o_{S_{\nb(t)}} = \left\{
  \begin{array}{ll}
    \ds \sum_{\ell=1}^{ S_{n(t)} - S_{\nb(t)} } \!\! \z_{S_{\nb(t)} + \ell} \,,
    & \mbox{if } S_{n(t)} > S_{\nb(t)} \,; \\[17pt]
    0 \,, & \mbox{if } S_{n(t)} = S_{\nb(t)} \,; \\[3pt]
    \ds \ - \hspace{-18pt} \sum_{\ell=0}^{ S_{n(t)} - S_{\nb(t)} + 1} \!\!\!\!\!
    \z_{S_{\nb(t)} + \ell} \,,
    & \mbox{if } S_{n(t)} < S_{\nb(t)} \,.
  \end{array}
  \right.
\ee
Therefore, using also (\ref{new-40}),
\be\label{new-60}
\begin{split}
  E_\o \! \left( 1_G \left| Y_{n(t)} - Y_{\nb(t)} \right| \right)
  & < E_\o \! \left( 1_G \!\! \sum_{\ell=0}^{S_{n(t)} - S_{\nb(t)}}
  \! \z_{S_{\nb(t)} + \ell} \right) \\
  & \le \left( C \sqrt{ \phi(t) } + 1 \right)  \sup_{|\ell| \le C \sqrt{ \phi(t) }
  + 1} \!\!\! E_\o \! \left( \z_{S_{\nb(t)} + \ell} \right) \,.
\end{split}
\ee
Since, for $t \to \infty$, $\nb(t) \sim t$ and $\phi(t) = o(t)$,
Lemma \ref{lem-mom} can be applied to the leftmost term of
(\ref{new-60}). Accordingly,
(\ref{new-55}) and (\ref{new-60}) imply
\be\label{new-70}
  E_\o \! \left( 1_G \left| f \! \left(  \frac{|Y_{n(t)} - Y_{\nb(t)} |}
  {\sqrt{t}} \right) -f(0) \right| \right) \le C' \sqrt{ \frac{\phi(t)} t }
  \le \frac{\e}3 \,,
\ee
for some constant $C' > 0$ and all $t$ large enough. This,
together with (\ref{new-50}), gives (\ref{new-20}), and concludes
the proof of Theorem \ref{Th:qCLT}.
\qed

\appendix

\section{Appendix: Technical lemmas}
\label{app-tech}

\subsection{Ergodicity of the PVP \dsy}
\label{app-pvp}

In this section we give the prove the \erg ity of the PVP
\dsy\ introduced in Section \ref{subs-pvp}, and another
related result.
\bigskip

\proofof{Theorem \ref{thm-pvp}} We follow the same ideas as in
\cite{L2, L3, CLS}. Let us first
prove that $T$ preserves $\nu$.

Set $A := B \times C$, where $B$ is an elementary cylinder of $\Z^{\Z^+}$
and $C$ is a measurable subset of $(\R^+)^\Z$. It is not hard to see that
$T^{-1}(A) = \bigsqcup_{k \in \Z} B_k \times \sigma_\z^{-k} (C)$,
where
\be\label{pvp-4}
B_k := \lset{(k, \x_1, \x_2, \ldots) \in \Z^{\Z^+}} {(\x_1, \x_2, \ldots) \in B} \, .
\ee
By the choice of $B$ and by definition of $Q_o$, $Q_o (B_k) =
p_k Q_o(B)$. Also, by the $P_o$-invariance of $\sigma_\z$,
$P_o (\sigma_\z^{-k} (C)) = P_o (C)$. This shows that
\be\label{pvp-6}
\nu(T^{-1}(A))  = \sum_{k \in \Z} \nu( B_k \times \sigma_\z^{-k} (C) )
= \sum_{k \in \Z} p_k \, Q_o(B) P_o (C) = \nu(A).
\ee
This then extends to all measurable sets $A$, proving our first
assertion. For the second assertion we need a lemma.

\begin{lemma}
\label{lem-pvp}
  Every $T$-invariant set $A \subseteq \Sigma$ is of the form
  $A = \Z^{\Z^+} \times C$ mod $\nu$ (meaning that the equality holds
  up to $\nu$-null sets),
  where $C$ is a measurable set of $(\R^+)^\Z$.
\end{lemma}

\proofof{Lemma \ref{lem-pvp}} We first give some preliminary definitions
and results. Let us endow $\Z^{\Z^+}$ with the distance
\be
  d(\x, \x') := \left[ \min\! \rset{n \in \Z^+} {\x_n \ne \x'_n} \right]^{-1}.
\ee
This is an \emph{ultrametric} distance, namely,
$\forall \x, \x', \x'' \in \Z^{\Z^+}$,
\be
  d(\x, \x'') \le \max \{ d(\x, \x') , d(\x', \x'') \} \,,
\ee
and its (open) balls are the cylinders
\be \label{pvp-10}
  \cB_\eps(\x) := \rset{\x'  \in \Z^{\Z^+}} {\x_i' = \x_i, \ \forall i = 1, 2,
  \ldots, [ \eps^{-1} ] } \, ,
\ee
where, once again, $[ \cdot ]$ indicates the integer part of a real number.
This makes $\Z^{\Z^+}$ a Polish ultrametric space, which is an
observation that will soon be useful.

Given $\cB_\eps(\x)$, as in (\ref{pvp-10}), and an elementary cylinder
$B \subseteq
\cB_\eps(\x)$, namely $B = \rset{\x'} {\x'_i = \x_i, \ \forall i = 1, 2, \ldots,
k}$, with $k \ge [ \eps_{-1} ]$, it is easy to see that
\be \label{pvp-20}
  \frac{ Q_o( \sigma_\x(B) ) } { Q_o( \sigma_\x(\cB_\eps(\x)) ) } =
  \frac{ Q_o(B) } { Q_o( \cB_\eps(\x) ) } \, .
\ee
So this holds for any measurable $B \subseteq \cB_\eps(\x)$.
If $B$ is not necessarily a subset of $\cB_\eps(\x)$, we can only
state that
\be \label{pvp-30}
  Q_o( \sigma_\x(B) \,|\, \sigma_\x(\cB_\eps(\x)) ) \ge
  Q_o(B \,|\, \cB_\eps(\x) ) \,.
\ee
(This follows from (\ref{pvp-20}), replacing $B$ with
$B \cap \cB_\eps(\x)$ and using the general inclusion
$\sigma_\x( B \cap \cB_\eps(\x) ) \subseteq \sigma_\x(B) \cap
\sigma_\x(\cB_\eps(\x))$.)

Now, recall the definitions of $F_\z$ and $Q_o (\cdot \,|\, F_\z)$
from the statement of Lemma \ref{lem-pvp-me}. It follows
from the above arguments that $F_\z$ is a Polish ultrametric
space endowed with the Borel \me\ $Q_o (\cdot \,|\, F_\z)$.
By \cite[Prop. 2.10]{M}, Lebesgue's Density Theorem holds.
In particular, if $\cB_\eps(\x,\z)$ denotes the ball of center $(\x,\z)$
and radius $\eps$ in $F_\z$ (corresponding to $\cB_\eps(\x)
\subset \Z^{\Z^+}$
through the identification $F_\z \cong \Z^{\Z^+}$), we have the
following:

\begin{lemma} \label{ldt}
  Let $A$ be a measurable subset of $\Sigma$ and $\z \in (\R^+)^\Z$
  be such that $A \cap F_\z$ is measurable (this happens for
  $P_o$-a.e.\ $\z$). Then, a.e.\ $(\x, \z) \in A \cap F_\z$, relative to
  $Q_o (\cdot \,|\, F_\z)$, is a density point of $A \cap F_\z$. This
  means that
  \begin{displaymath}
    \lim_{\e \to 0^+} Q_o ( A \,|\, \cB_\eps(\x, \z) ) :=
    \lim_{\e \to 0^+} \frac{ Q_o ( A \cap \cB_\eps(\x, \z) \,|\, F_\z) }
    { Q_o ( \cB_\eps(\x, \z) \,|\, F_\z) } = 1 \,.
  \end{displaymath}
\end{lemma}

We finally come to the actual proof of the lemma. Let us first assume
$\nu(A) > 0$, otherwise one sets $C := \emptyset$ and
the proof is finished. Then, by contradiction, we assume that $A$ is not of
the type $A = \Z^{\Z^+} \times C$ mod $\nu$, that is, it is not a union
of horizontal fibers, modulo null sets. Therefore, for a small
enough $\delta>0$, the set
\be \label{pvp-35}
  C_\delta := \lset{\z \in (\R^+)^\Z} {\delta \le Q_o (A \,|\, F_\z)
  \le 1-\delta}
\ee
has positive $P_o$-\me. Set $A' := \bigsqcup_{\z \in C_\delta}
(A \cap F_\z)$. By Fubini, $\nu(A')>0$.

We claim that we can find a point $(\x, \z) \in A'$ which is
both a recurrent
point to $A'$ w.r.t.\ $T$ (i.e.,
\be \label{pvp-40}
  T^n (\x, \z) \in A' \,,
\ee
for countably many values of $n$), and a density point
of $A \cap F_\z$, relative to $Q_o ( \, \cdot \,|\, F_\z)$.
This is true because, by
Poincar\'e's Recurrence Theorem, $\nu$-a.a.\ points in $A'$
recur to $A'$. By Fubini and the definition of $A'$, this implies that,
for $P_o$-a.a.\
$\z \in C_\delta$, $(\x, \z)$ is a recurrent point (to $A'$), for $Q_o$-a.e.\
$\x \in A' \cap F_\z$. Now, consider a typical $\z$ in the sense just
described and exclude from the recurrent points contained in
$A' \cap F_\z$ those
that are not density points of $A' \cap F_\z$. By Lemma
\ref{ldt} this amounts
to excluding a negligible set of points. Any remaining point
verifies our claim (in fact, $A' \cap F_\z = A \cap F_\z$,
by definition of $A'$).

Therefore, we can find a large enough $n$ that
verifies (\ref{pvp-40}) and
\be \label{pvp-50}
  Q_o(A  \,|\, \cB_{1/n} (\x,\z) ) > 1 - \delta \, .
\ee
Notice that, via (\ref{T}) and (\ref{underRW}), it is easy to find
an expression for the iterates of $(\x,\z)$:
\be \label{Tn}
  T^n(\x, \z) := (\sigma_\x^n (\x), \sigma_\z^{S_n} (\z)) \,.
\ee
The above makes it clear that $T^n$ acts on $F_\z$ by operating
$n$ shifts in the dynamical variable and mapping the environment
to the new environment $\sigma_\z^{S_n} (\z)$. But, by
(\ref{pvp-10}), $\sigma_\x^n (\cB_{1/n} (\x)) = \Z^{\Z^+}$. Therefore,
$T^n (\cB_{1/n} (\x,\z)) = F_{\sigma_\z^{S_n} (\z)}$.
On the other hand, using the invariance of $A$, (\ref{pvp-30}) and
(\ref{pvp-50}), we can write
\be
\begin{split}
  Q_o ( A \,|\, F_{\sigma_\z^{S_n} (\z)} ) &=
  Q_o ( T^n(A) \,|\, T^n (\cB_{1/n} (\x,\z)) ) \\
  &\ge Q_o ( A \,|\, \cB_{1/n} (\x,\z) ) \\
  &> 1 - \delta \,,
\end{split}
\ee
which, in view of (\ref{pvp-35}), shows that $\sigma_\z^{S_n} (\z)
\not\in C_\delta$. But (\ref{pvp-40}) and (\ref{Tn}) imply that
$\sigma_\z^{S_n} (\z) \in C_\delta$, which is the sought
contradiction.

Therefore $A = \Z^{\Z^+} \times C$ mod $\nu$, for some $C \subseteq
(\R^+)^\Z$. We still need to
prove that $C$ is measurable. If not, by \cite[Lem. A.1]{L1}, there
exists a measurable $C'$ such that $C \symmdiff C'$ is contained in
a null set, implying $A = \Z^{\Z^+} \times C'$ mod $\nu$.
\qed

To end the proof of Theorem \ref{thm-pvp} suppose,
again by contradiction,
that the \sy\ has an invariant set $A$, which, by the above lemma,
must be of the form $A = \Z^{\Z^+} \times C$ mod $\nu$,
with $0 < P_o(C) < 1$. By the \erg ity of
$\sigma_\z$, there must be a subset $C' \subseteq C$, with
\be \label{pvp-60}
P_o(C') > 0 \,,
\ee
such that $\sigma_\z (C') \subseteq C^c := (\R^+)^\Z \setminus C$.

For $\z \in C'$, set $B_{\z,1} := \rset{(\x,\z) \in F_\z}
{\x_1 = 1}$. Then
\be \label{pvp-70}
Q_o \left( B_{\z,1} \,|\, F_\z \right) = Q_o\! \left(
\lset{\x \in \Z^{\Z^+}} {\x_1 =1} \right) = p_1 > 0 \,,
\ee
by the assumptions on $(p_k)$ (symmetry, half-monotonocity,
and positive variance; see, however, Remark \ref{rk-pvp}).
Also, $T (B_{\z,1}) = F_{\sigma_\z(\z)}$ and $\sigma_\z(\z) \in C^c$.
Therefore, setting $A_o := \bigsqcup_{\z \in C'} B_{\z,1}$,
one has that $T(A_o) \subseteq \Z^{\Z^+} \times C^c$,
with $\nu(A_o) > 0$ (the latter inequality coming from
(\ref{pvp-60}), (\ref{pvp-70}) and Fubini's Theorem).
This contradicts the invariance of $A$
and thus proves the theorem.
\qed

\skippar

\proofof{Lemma \ref{lem-pvp-me}}
The proof of Theorem \ref{thm-pvp} (see in particular
(\ref{pvp-4})-(\ref{pvp-6}) and the concluding paragraph) shows
that $T$ maps $B_{\z,k} := \rset{(\x,\z) \in F_\z}
{\x_1 = k}$ onto $F_{\sigma_\z^k (\z)}$, pushing the \me\
$Q_o ( \,\cdot\, | F_\z )$, restricted to $B_{\z,k}$, to $p_k
Q_o ( \,\cdot\, |  F_{\sigma_\z^k (\z)} )$. Since $F_\z =
\bigsqcup_{k \in \Z} B_{\z,k}$, both statements of the lemma
follow.
\qed

\begin{remark}
  \label{rk-pvp}
  The proof of Lemma \ref{lem-pvp-me} helps
  to show that Theorem \ref{thm-pvp}
  holds under much weaker assumptions on the underlying \rw:
  it suffices to require that $v_p > 0$. This, in fact, implies that
  $p_k > 0$, for some $k \ne 0$. The proof of Theorem \ref{thm-pvp}
  still \fn s if, in the last two paragraphs, one substitutes $\sigma_\z$
  with $\sigma_\z^k$ (also \erg) and $B_{\z,1}$ with $B_{\z,k}$.
\end{remark}

\subsection{Averaging}
\label{app-avg}

The next lemma, which is needed in the proof of the main theorems (cf.\
Lemma \ref{lem-mom}), proves an assertion that can be roughly
described as follows: given a \fn\ $a: \Z \into \R$ and an ``expanding''
sequence of \pr\ densities on $\Z$ that are increasing on $\Z^-$ and
decreasing on $\N$, the expected value of $a$ relative to these
densities tends to its Cesaro average.

\begin{lemma}
  \label{lem-tech1}
  Let $a := (a_j, \, j \in \Z) \subset \R$ be such that
  \begin{displaymath}
    \lim_{k \to \infty} \, \frac1 k  \sum_{j=0}^{k-1} a_j =
    \lim_{k \to \infty} \, \frac1 k \sum_{j=-1}^{-k} a_j =
    \bar{a} \,.
  \end{displaymath}
  For $n \in \N$, let $p^{(n)} = (p^{(n)}_j, \, j \in \Z)$ be the density of
  a \pr\ distribution on $\Z$, whose expectation we denote by
  $\mathcal{E}_n$, such that:
  \begin{itemize}
  \item[(i)] $j \mapsto p^{(n)}_j$ is increasing on $\Z^-$ and decreasing
  on $\N$;
  \item[(ii)] for all $r \in \N$, $\mathcal{E}_n ( 1_{[-r,r]} ) :=
  \sum_{j=-r}^r p^{(n)}_j$ vanishes as $n \to \infty$.
  \end{itemize}
  Then
  \begin{displaymath}
    \lim_{n \to \infty} \mathcal{E}_n (a) := \lim_{n \to \infty}
    \sum_{j \in \Z} p^{(n)}_j a_j = \bar{a} \,.
  \end{displaymath}
\end{lemma}

\proof For $k_1, k_2 \in \Z$, $k_1 \le k_2$, set
\begin{equation}
  \mathcal{U}_{k_1,k_2} (a) := \frac1 {k_2 - k_1 + 1} \,
  \sum_{j=k_1}^{k_2} a_j \,.
\end{equation}
Define also
\begin{align}
  \label{tech1-20}
  \mathcal{E}_n^+ (a) := & \sum_{j=0}^\infty p^{(n)}_j a_j \,, \\
  \label{tech1-30}
  \mathcal{E}_n^- (a) := & \sum_{j=-1}^{-\infty} p^{(n)}_j a_j \,.
\end{align}
Let us approximate (\ref{tech1-20}) and (\ref{tech1-30}) separately.
It is not hard to see (by ``slicing'' the density $p^{(n)}$ horizontally)
that
\begin{equation}
  \label{tech1-40}
  \mathcal{E}_n^+ (a) = \sum_{j=0}^\infty \left( p^{(n)}_j - p^{(n)}_{j+1}
  \right) (j + 1) \, \mathcal{U}_{0,j} (a) \,.
\end{equation}

Fix $\eps>0$. The hypotheses on $a$ imply that
$\exists r \in \N$ so large that, $\forall j > r$,
\begin{align}
  \label{tech1-50}
  | \mathcal{U}_{0,j} (a) - \bar{a} | & \le \eps/2 \,; \\
  \label{tech1-60}
  | \mathcal{U}_{-j,-1} (a) - \bar{a} | & \le \eps/2 \,.
\end{align}
Set
\begin{equation}
  \mathcal{E}_{n,r}^+ (a) := \sum_{j=r+1}^\infty \left( p^{(n)}_j -
  p^{(n)}_{j+1} \right) (j + 1) \mathcal{U}_{0,j} (a) \,.
\end{equation}
(\ref{tech1-50}) implies that
\begin{equation}
  \label{tech1-80}
  \left| \mathcal{E}_{n,r}^+ (a) - \bar{a} \, \mathcal{E}_{n,r}^+ (1) \right|
  \le \frac{\eps}2 \, \mathcal{E}_{n,r}^+ (1) \,,
\end{equation}
with the understandable meaning that $1$ also denotes the
sequence that is identically equal to 1.

Analogously, the term (\ref{tech1-30}) can be rewritten as
\begin{equation}
  \mathcal{E}_n^- (a) = \sum_{j=1}^\infty \left( p^{(n)}_{-j} - p^{(n)}_{-j-1}
  \right) j \, \mathcal{U}_{-j,-1} (a) \,,
\end{equation}
and, upon defining
\begin{equation}
  \mathcal{E}_{n,r}^- (a) := \sum_{j=r+1}^\infty \left (p^{(n)}_{-j} -
  p^{(n)}_{-j-1} \right) j \, \mathcal{U}_{-j,-1} (a),
\end{equation}
we get through (\ref{tech1-60}) that
\begin{equation}
  \label{tech1-110}
  \left| \mathcal{E}_{n,r}^- (a) - \bar{a} \, \mathcal{E}_{n,r}^- (1) \right| \le
  \frac{\eps}2 \, \mathcal{E}_{n,r}^- (1) \,.
\end{equation}
If we name $\mathcal{E}_{n,r} (\cdot) := \mathcal{E}_{n,r}^+ (\cdot) +
\mathcal{E}_{n,r}^- (\cdot)$, we obtain from (\ref{tech1-80}) and
(\ref{tech1-110}) that
\begin{equation}
  \label{tech1-120}
  \left| \mathcal{E}_{n,r} (a) - \bar{a} \, \mathcal{E}_{n,r} (1) \right| \le
  \frac{\eps}2 \, \mathcal{E}_{n,r} (1) \le \frac{\eps}2 \,.
\end{equation}

On the other hand, it is clear from the above arguments that
$1 - \mathcal{E}_{n,r} (1) = \mathcal{E}_n (1) - \mathcal{E}_{n,r} (1)$
is a portion of the mass of $p^{(n)}$ contained in $[-r, r]$, which is
measured by $\mathcal{E}_n (1_{[-r, r]})$. Therefore, defining
$\cA := \max_{|j| \le r} |a_j|$ and using \emph{(ii)}, there exists $N =
N(\eps, \bar{a}, \cA, r)$ such that, for all $n \ge N$,
\begin{equation}
  \label{tech1-130}
  1 - \mathcal{E}_{n,r} (1) \le \mathcal{E}_n (1_{[-r, r]} ) \le
  \frac{\eps} {2 (\cA + |\bar{a}| )} \,.
\end{equation}
Notice that $N$ can be thought of as a \fn\ of $\eps$ and the sequence
$a$ (for $r$ is also a \fn\ of $\eps$ and $a$). Finally, $\forall n \ge N$,
\begin{equation}
\begin{split}
 \left| \mathcal{E}_n (a) - \bar{a} \right| &\le
     \left| \mathcal{E}_n (a) - \mathcal{E}_{n,r} (a) \right|
  + \left| \mathcal{E}_{n,r} (a) - \bar{a} \, \mathcal{E}_{n,r} (1) \right|
  + \left| \bar{a} \, \mathcal{E}_{n,r} (1) - \bar{a} \right| \\
  &\le ( \cA + |\bar{a}| ) \left( 1 - \mathcal{E}_{n,r} (1) \right)
  + \left| \mathcal{E}_{n,r} (a) - \bar{a} \, \mathcal{E}_{n,r} (1) \right| \\
  &\le \eps/2 + \eps/2 \,,
\end{split}
\end{equation}
by (\ref{tech1-120}) and (\ref{tech1-130}). This completes the proof.
\qed

In the main body of the paper, Lemma \ref{lem-tech1} is used
with $p^{(n)}$ being the \pr\ density of the underlying
\rw\ at time $n$. In order to show that such densities verify condition
\emph{(i)} above, we need another simple lemma.

\begin{lemma}
  \label{lem-tech2}
  If $p$ and $p'$ are symmetric and half-monotonic densities on
  $\Z$ (see definitions in Section \ref{subs-model}), so is their
  convolution $p * p'$.
\end{lemma}

\proof Let us first treat the special case $p = q^{(r)}$ and $
p' = q^{(r')}$, with $r, r' \in \N$, where
\begin{equation}
  \label{tech2-20}
  q^{(r)}_k = \left\{
  \begin{array}{ll}
    (2r+1)^{-1} \,, & -r \le k \le r \,; \\
    0 \,, & |k| > r \,.
  \end{array} \right.
\end{equation}
Assume $r \ge r'$: this is no loss of generality as the convolution
is symmetric. It is easy to calculate that
\begin{equation}
  \label{tech2-30}
  \left( q^{(r)} * q^{(r')} \right)_j = \frac1 {2r +1} \cdot
  \left\{
  \begin{array}{ll}
    1 \,, & |j| \le r-r' \,; \\[3pt]
    \ds \frac{r+r'+1-j} {2r'+1} \,, & r-r' < |j| \le r+r'
    \,; \\[3pt]
    0 \,, & |j| > r+r' \,,
  \end{array} \right.
\end{equation}
which is symmetric and half-monotonic.
Now, a general $p$ as in the statement of the lemma can be rewritten
as $p = \sum_{r=0}^\infty (p_r - p_{r+1}) (2r + 1) q^{(r)}$, and
analogously for $p'$. Hence
\begin{equation}
  \label{tech2-40}
  p * p' = \sum_{r=0}^\infty \sum_{r'=0}^\infty
  (p_{r} - p_{r+1}) (p'_{r'} - p'_{r'+1}) (2r + 1) (2r' + 1) \,
  q^{(r)} * q^{(r')} ,
\end{equation}
which is symmetric and half-monotonic because it is a
countable linear combination
of symmetric and half-monotonic densities,
with positive coefficients.
\qed

\subsection*{Compliance with ethical standards}

\smallskip\noindent
{\bf Acknowledgments and funding.}
We would like to thank two referees for carefully reading the
first version of the manuscript and making several useful comments.
Work in Bologna has been partially supported by the FIRB
Projects RBFR08UH60 and RBFR10N90W (MIUR, Italy). A.~Bianchi
acknowledges Universit\`a di Padova for partial support
through the Project \emph{``Stochastic Processes and Applications to
Complex Systems''} (CPDA123182). M.~Ligab\`o acknowledges the
Fondazione Cassa di Risparmio di Puglia and the Gruppo Nazionale di
Fisica Matematica (INdAM, Italy) for partial support. G.~Cristadoro
and M.~Lenci also declare that this work is part of their activity within
the Gruppo Nazionale di Fisica Matematica.

\medskip\noindent
{\bf Conflict of interest.} The authors declare that they have no conflict
of interest.

\end{document}